\documentclass[a4paper,10pt]{amsart}
\usepackage{amssymb,amsfonts,amsmath,amsthm}
\usepackage{enumerate}
\usepackage{epic,eepic}
\usepackage{graphicx}

\advance \oddsidemargin by -1in
\advance \evensidemargin by -1in
\advance \topmargin by -1in
\newtheorem{thm}{Theorem}
\newtheorem{lem}{Lemma}
\newtheorem{cor}{Corollary}
\theoremstyle{remark}
\newtheorem{exmp}{Example}
\newtheorem{defn}{Definition}

\def\R{\mathbb{R}}
\def\N{\mathbb{N}}
\def\C{\mathbb{C}}
\def\Z{\mathbb{Z}}
\def\cd{\mathcal{D}}
\def\cf{\mathcal{F}}
\def\cj{\mathcal{J}}
\def\cc{\mathcal{C}}
\def\co{\mathcal{O}}

\def\VG{V \mkern -2.4mu G}
\def\VhatC{V \mkern -2.4mu \hat C}
\def\EG{E \mkern -0.2mu G}
\def\EhatC{E \mkern -0.2mu \hat C}
\def\VX{V \mkern -2.9mu X}
\def\VH{V \mkern -2.9mu H}
\def\VK{V \mkern -2.9mu K}
\def\EX{E \mkern -2mu X}
\def\VY{VY}
\def\VtildeY{V \tilde Y}
\def\EY{EY}
\def\VhatD{V \mkern -3mu \hat D}

\DeclareMathOperator{\Aut}{Aut}
\DeclareMathOperator{\diam}{diam}
\def\Im{\qopname\relax o{Im}}

\def\eps{\varepsilon}
\def\abs#1{\lvert #1\rvert}

\newcounter{fig}
\begin{document}

\title[Transition probabilities on self-similar graphs]
{Asymptotics of the transition probabilities of the simple random walk on self-similar graphs}


\author[B.~Kr\"on]{Bernhard~Kr\"on$^1$}
\address{Bernhard Kr\"on\\
Erwin Schr\"odinger Institute (ESI) Vienna\\
Boltzmanngasse 9\\
1090 Wien\\
Austria}
\email{bernhard.kroen@univie.ac.at}

\author[E.~Teufl]{Elmar~Teufl$^2$}
\address{Elmar Teufl\\
Department of Mathematics C\\
Graz University of Technology\\
Steyrergasse 30\\
8010 Graz\\
Austria}
\email{elmar.teufl@tugraz.at}

\thanks{$^1$Bernhard~Kr\"on is supported by the project P14379-MAT of the Austrian Science Fund (FWF)}
\thanks{$^2$Elmar~Teufl is supported by the START-project Y96-MAT of the FWF}

\begin{abstract}
It is shown explicitly how self-similar graphs
can be obtained as `blow-up' constructions of finite cell graphs $\hat C$.
This yields a larger family of graphs than the graphs obtained by discretising
continuous self-similar fractals.

For a class of symmetrically self-similar graphs
we study the simple random walk on a cell graph $\hat C$,
starting in a vertex $v$ of the boundary of $\hat C$.
It is proved that the expected number of returns to $v$
before hitting another vertex in the boundary coincides
with the resistance scaling factor.

Using techniques from complex rational iteration and singularity analysis
for Green functions we compute the asymptotic behaviour
of the $n$-step transition probabilities of the simple random walk on the whole graph.
The results of Grabner and Woess for the Sierpi\'nski graph are generalised to
the class of symmetrically self-similar graphs and
at the same time the error term of the asymptotic expression is improved.
Finally we present a criterion for the occurrence of oscillating phenomena
of the $n$-step transition probabilities.
\end{abstract}

\date{\today}

\maketitle

\section{Introduction}\label{sec:intro}

Self-similar fractals are usually constructed as the compact invariant set
of an iterated function system.
Analysis on these sets is still a rapidly growing field in mathematics.
One approach is to study self-similar graphs as the discretization of fractals
and then transfer the results on these graphs back to the fractals via rescaling.
In this way one can construct Brownian motion on fractals by starting
with random walks on graphs, see for example Barlow \cite{barlow98diffusions}, Lindstr{\o}m
\cite{lindstroem90brownian} and the references therein.

Contracted images of a fractal as compact subsets
of the fractal are called \emph{cells}.
Correspondingly, self-similar graphs have \emph{cell graphs} as finite subgraphs
which carry essential information of the whole graph.
Malozemov and Teplyaev defined self-similarity
of graphs axiomatically in \cite{malozemov95pure}.
Their definition was restricted to the case where the cells of the graph have
exactly two boundary points. In \cite{kroen02green} one of the authors introduced
self-similarity for an arbitrary number of boundary points.
Another axiomatic approach for self-similar graphs with more boundary points
was chosen in \cite{malozemov01self}.
Sabot constructed self-similar graphs using equivalence relations on word spaces,
see \cite{sabot01spectral}.

Probably the most extensively studied self-similar graph is the Sierpi\'nski graph.
In \cite{grabner97functional1} Grabner and Woess considered the Green function
of the simple random walk on the Sierpi\'nski graph which describes the returning
of the random walk to a certain origin vertex.
A combinatorial substitution based on path arguments was used
to obtain a functional equation for this Green function.
Iterating this equation they obtained the analytic continuation
of the Green function as a rapidly converging product of rational terms.
The asymptotic behaviour and oscillation phenomena were computed
by analysing the singularity of the Green function in $z=1$.
These techniques had been introduced by Flajolet and Odlyzko
in \cite{flajolet90singularity,odlyzko82periodic}.
\emph{Symmetrically self-similar} graphs were constructed in \cite{kroen02green}
as class of self-similar graphs such that the combinatorial
substitution mentioned above can be applied.
Sections~\ref{sec:asymptotic} to \ref{sec:transition} of the present paper
are a generalisation and a further development
of the asymptotic analysis in \cite{grabner97functional1}.

In Section~\ref{sec:ssg} we recall from \cite{kroen02green} the axiomatic definition
of self-similarity of graphs and a graph theoretic analogue
to the Banach fixed point theorem for self-similar graphs.
Heuristically, this theorem says that a self-similar graph has a ``centre'',
either as an origin vertex or as an origin cell.
We will later focus on self-similar graphs with an origin vertex.
Our axiomatic approach to self-similarity of graphs is based
on contracting a given infinite graph.
On the other hand one can obtain self-similar graphs by a ``blow-up'' construction,
starting with a finite cell graph.
Theorem~\ref{thm:blowup} yields this construction explicitly,
the resulting graphs may have more irregular structures
than discretized post critically finite self-similar (pcfss) sets,
see \cite{barlow98diffusions,kigami93harmonic} for the notion of pcfss sets.
It is proved that a symmetrically self-similar graph is bipartite if and only if
its cell graphs are bipartite.

Section~\ref{sec:green} is devoted to the study of Green functions.
For the Sierpi\'nski graph, a functional equation
\[ G(z) = f(z) G(d(z)) \]
is valid  for the Green function $G$ at the origin vertex, which
was first observed by Rammal and Toulouse \cite{rammal84random,rammal83random},
see also \cite{grabner97functional1}.
In \cite{kroen02green} this equation was generalised to the class of
so called \emph{symmetrically} self-similar graphs.
One property of these graphs is that their cell graphs are all isomorphic.
The functions $d$ and $f$ are rational functions associated with
the simple random walk on the cell graph $\hat C$.
The \emph{transition function} $d$ is the generating function
of the probabilities that the simple random walk,
starting in the boundary of $\hat C$, hits another vertex
in the boundary for the first time after exactly $n$ steps. The product
\[ \prod_{k=0}^{\infty} f(d^k(z)) \]
converges on the Fatou set $\cf$ of $d$
(possibly with the exception of countably many points)
to the Green function $G$ as a solution of the functional equation above.
The point $z=0$ is the only attracting fixed point of $d$ and
the Fatou set $\cf$ contains $(\C \setminus \R) \cup (-1,1)$.

The simple random walk on a cell graph $\hat C$
of a symmetrically self-similar graph $X$
is studied in Section~\ref{sec:geometric}.
The number $\theta$ of vertices in the boundary of a cell graph
does not depend on the choice of the cell.
The graph $\hat C$ consists of $\mu$ amalgamated copies
of the $\theta$-complete graph $K_\theta$.
This parameter $\mu$ corresponds to the usual
\emph{mass scaling factor} of self-similar sets.
The number $\tau=d'(1)$ is called \emph{time scaling factor}.
Let this random walk start at a vertex $v$ in the boundary of $C$.
Then $\tau$ is the average time until hitting another boundary vertex than $v$.
For symmetrically self-similar graphs with bounded geometry
it is proved that $f(1) = \rho$, where $f(1)$
is the expected number of returns to $v$
before hitting a new vertex in the boundary of $\hat C$
and $\rho$ is the \emph{resistance scaling factor}, see \cite{barlow98diffusions}.
For the proof we use the relation $\tau = \mu\rho$.

In Section~\ref{sec:asymptotic} we study the asymptotic behaviour
of the $n$-step return probabilities.
Using the functional equation of Section~\ref{sec:green}
we derive a local singular expansion of $G$ in $z=1$,
\[ G(z) = (1-z)^{\frac{\log(\mu)}{\log(\tau)}-1}
        \biggl(\omega\biggl(\frac{\log(1-z)}{\log(\tau)}\biggr) + 
	\co_\delta(\abs{z-1})\biggr) \]
for all $z\in\C$ such that $\abs{z-1}<r$ and $\abs{\arg(1-z)}<\delta$,
where $0 < \delta < \pi$ and $r > 0$.
The function $\omega$ is $1$-periodic and holomorphic
on some horizontal strip around the real axes.
In \cite{grabner97functional1} this type of expansion
was computed for the Sierpi\'nski graph with a weaker error term.
Using the method of singularity analysis we obtain the asymptotics
of the $n$-step return probabilities $p^{(n)}(o,o)$ to an origin vertex $o$.
For the case that $X$ is not bipartite we have
\[ p^{(n)}(o,o) = n^{-\frac{\log(\mu)}{\log(\tau)}}
        \biggl( \sigma\biggl(\frac{\log(n)}{\log(\tau)}\biggr) +
	\co(n^{-1}) \biggr), \]
where $\sigma$ is a $1$-periodic, holomorphic function defined
on some horizontal strip around the real axes.
If $X$ is bipartite, then $p^{(2n+1)}(o,o)=0$ and
$p^{(2n)}(o,o)$ satisfies the above asymptotic behaviour.
This type of oscillation seems to be typical
for random walks on self-similar graphs and Brownian motion on fractals,
see \cite{barlow88brownian,arous00large,fukushima92spectral,
grabner97functional2,grabner97functional1,kigami93weyl,teufl01average}.


It is proved in Section~\ref{sec:oscillation}
that the function $\sigma$ is not constant
if the Julia set $\cj$ of $d$ is a Cantor set.

In Section~\ref{sec:transition} it is shown that
for any pair of vertices $x$ and $y$, the transition probabilities $p^{(n)}(x,y)$
satisfy the same asymptotic behaviour as $p^{(n)}(o,o)$.
For this purpose we use a ratio limit theorem,
see \cite{grabner97functional1}.

Section \ref{sec:line} is devoted to the case
where the cell graph is a finite line $L_n$ of length $n$.
The corresponding transition functions $d_n$ are conjugated
to the Chebychev polynomials $T_n$ via conjugacy map $z \mapsto \frac{1}{z}$.
Their Julia set is $(\R \cup \{\infty\}) \setminus (-1,1)$
and the periodic function $\sigma$ is constant.
We study the probability for hitting a given vertex
in the minimal number of steps for general, locally finite, reversible Markov chains.
This probability is maximal whenever the process can be projected
on the simple random walk on $L_n$.
From this we deduce that the cell-graphs of symmetrically self-similar graphs
are isomorphic to $L_n$ if and only if
the transition function is conjugated to $T_n$ in the sense above.
We are interested in characterisations of this type
because we believe that these are the only symmetrically self-similar graphs
with a constant periodic function $\sigma$.
This would mean that in all other cases we may observe a non trivial oscillation
of the asymptotic transition probabilities.

In Section~\ref{sec:examples} some examples are discussed.

\section{Self-similar graphs}\label{sec:ssg}

Graphs $X=(\VX,\EX)$ with vertex set $\VX$ and edge set $\EX$ are always
supposed to be undirected, without loops or multiple edges.
Vertices $x$ and $y$ are \emph{adjacent} if $\{x,y\}$ is an edge in $\EX$.
The \emph{degree} $\deg_X (x)$ of $x$ is the number of vertices in $\VX$ being adjacent to $x$.
A graph has \emph{bounded geometry} if the set of vertex degrees is bounded.
The \emph{vertex boundary} or \emph{boundary} $\theta C$ of a set $C$ of vertices in $\VX$
is the set of vertices in $\VX \setminus C$ which are adjacent to some vertex in $C$.
A \emph{path of length $n \in \N_0$} from a vertex $v$ to a vertex $w$ is
an $(n+1)$-tuple $\pi=(v=x_0,x_1,\ldots,x_n=w)$ of vertices
such that $x_{i-1}$ and $x_i$ are adjacent for $1 \leq i \leq n$.
The path $\pi$ is \emph{closed} if $v=w$. The distance $d_X(v,w)$ of $v$ and $w$ in $X$ is
the length of the shortest path from $v$ to $w$.
A set $C$ of vertices is called \emph{connected} if any pair of vertices in $C$
can be connected by a path in $X$ which does not leave $C$.
We write $\hat C$ for the subgraph of $X$ spanned by $C \cup \theta C$.

We briefly repeat the definition of self-similar graphs,
see Definitions~1 and 2 in \cite{kroen02green} or in \cite{kroen02growth}:
Let $F$ be a set of vertices in $\VX$ and
let $\cc_X(F)$ be the set of connected components in $\VX \setminus F$.
We define the \emph{reduced graph} $X_F$ of $X$ be setting $\VX_F = F$
and connecting two vertices $x$ and $y$ in $\VX_F$ by an edge
if and only if there exists a component $C$ in $\cc_X(F)$
such that $x$ and $y$ are in the boundary of $C$.

\begin{defn}\label{def:selfsimilar}
A connected infinite graph $X$ is called \emph{self-similar}
with respect to $F \subset \VX$ and $\psi : \VX \to \VX_F$ if
\begin{enumerate}[(F1)]
	\item no vertices in $F$ are adjacent in $X$,
	\item the intersection of the boundaries of two different components in $\cc_X(F)$
		contains not more than one vertex and
	\item $\psi$ is an isomorphism between $X$ and $X_F$.
\end{enumerate}
We will write $\phi$ instead of $\psi^{-1}$
and $F^n$ instead of $\psi^n(F)$, where $\psi^n$ denotes the $n$-fold iterate of $\psi$.
Components of $\cc_X(F^n)$ are \emph{$n$-cells},
1-cells are also just called \emph{cells}.
An \emph{origin cell} is a cell $C$
such that $\phi(\theta C) \subset C$.
A fixed point of $\psi$ is called \emph{origin vertex}.
See Section~\ref{sec:examples} for examples.
\end{defn}

The following lemma is a graph theoretic analogue to the Banach fixed point theorem,
see \cite[Theorem~1]{kroen02green} or \cite[Lemma~1]{kroen02growth}.

\begin{lem}\label{lem:banach}
Let $X$ be self-similar with respect to $F$ and $\psi$.
Then $X$ is self-similar with respect to $F^k$ and $\psi^k$,
for any positive integer $k$. Either $X$ has
\begin{enum}
	\item exactly one origin cell and no origin vertex or
	\item exactly one origin vertex $o$.
		If in this case $X$ is locally finite,
		then there is a positive integer $n$ such that the subgraphs $X_A$ of $X$
		which are spanned by $A \cup \{o\}$, for components $A$ in $\cc_X(\{o\})$,
		are self-similar graphs with respect to
		\[ F_A = \{o\} \cup (F^n \cap A)
			\qquad\text{and}\qquad \psi_A = \psi^n|_{F_A}. \]
		The graphs $X_A$ have exactly one origin cell and $o$ is their origin vertex.
\end{enum}
\end{lem}

\begin{defn}
Let $X$ be a graph which is self-similar with respect to $F$ and $psi$.
We call $\hat C_n$ the \emph{$n$-cell graph} $C_n$,
1-cell graphs are also called \emph{cell graphs}.
The graph $X$ is \emph{symmetric} (or \emph{doubly symmetric})
if it is locally finite and
satisfies the following axioms:
\begin{enumerate}[(S1)]
	\item All cells are finite and for any pair of cells $C$ and $D$ in $\cc_X(F)$
		there exists an isomorphism $\alpha : \VhatC \to \VhatD$
		of $\hat C$ and $\hat D$ such that $\alpha(\theta C) = \theta D$.
	\item The automorphism group $\Aut(\hat C)$ of $\hat C$
		acts doubly transitive on $\theta C$,
		which means that it acts transitively on the set
		\[ \left\{ (x,y) \mid x,y \in \theta C, \; x \neq y \right\}, \]
		where $g((x,y))$ is defined as $(g(x),g(y))$ for any $g\in \Aut(\hat C)$.
\end{enumerate}
If a self-similar graph $X$ satisfies (S1) then
let $\theta$ be the number of vertices in the boundary of a cell.
Axiom (S1) implies that $\theta$ is independent of the choice of the cell $C$.
\end{defn}

The next theorem describes a `blow-up' construction
for self-similar graphs satisfying Axiom (S1).

\begin{thm}\label{thm:blowup}
Let $G=(\VG,\EG)$ be a finite, connected graph and
let $B$ be a subset of $\VG$ with the following properties:
\begin{itemize}
	\item The set $B$ has at least two elements, $\VG \setminus B$ is connected
		and no pair of vertices in $B$ is adjacent in $G$.
	\item The graph $G$ consists of $\mu$ complete subgraphs $K_\theta^1,\dotsc, K_\theta^\mu$
		which are isomorphic to the complete graph $K_\theta$ with $\theta$ vertices,
		such that $\theta=\abs{B}$ and any two of these complete subgraphs have at most
		one vertex in common.
\end{itemize}
Then for each $i$ in $\{1,\ldots,\mu\}$ there is a graph $X$,
a set $F\subset \VX$ and a map $\psi: \VX \to \VX_F$
such that $X$ is self-similar with respect to $F$ and $\psi$,
and the following statements are true:
\begin{enum}
	\item The set $\VG \setminus B$ is the unique origin cell of $X$,
		$G$ is the corresponding cell graph and $\psi(\VK_\theta^i)=B$.
	\item For any cell $C$ there exists an isomorphism $\alpha : \VhatC \to \VG$
		of $\hat C$ and $G$ such that $\alpha(\theta C) = B$.
	\item A vertex $o$ is the origin vertex of $X$
		if and only if it is a fixed point of $\psi$ in $\VK_\theta^i$.
		In this case $\VK_\theta^i \cap B = \{o\}$
		and $\VX \setminus \{o\}$ is connected.
\end{enum}
\end{thm}

\begin{proof}
We construct sequences of graphs $(Y_n)_{n \in \N}$,
sets of vertices $(F_n)_{n \in \N}$, where $F_n \subset \VY_n$,
and functions $(\psi_n)_{n \in \N}$ inductively such that they satisfy the following:
\begin{itemize}
\item The function $\psi_n$ is an isomorphism between $Y_{n-1}$
	and the reduced graph $(Y_n)_{F_n}$.
\item The graph $X = \bigcup_{n=1}^\infty Y_n$
	(seen as union of graphs and not as set theoretic union)
	is self-similar with respect to $F = \bigcup_{n=1}^\infty F_n$ and $\psi : \VX \to F$,
	where $\psi$ is defined by $\psi|_{\VY_n}=\psi_{n+1}$,
\item and $X$ has the properties required in the theorem.
\end{itemize}

First we set $\VY_{-1}=\emptyset$, $Y_0=K_\theta^i$, $Y_1=G$ and $F_1=B$.
Let $\psi_1 : \VY_0 \to F_1$ be any bijective function.
We suppose that there are graphs $Y_k$ and sets of vertices $F_k \subset \VY_k$,
for $1 \leq k \leq n$, such that $Y_{k-1}$ is a subgraph of $Y_k$ and
there is an isomorphism $\psi_k : \VY_{k-1} \to F_k$ of $Y_{k-1}$ and $(Y_k)_{F_k}$.

Let $\tilde Y_n$ be a disjoint, isomorphic copy of $Y_n$
and let $\alpha_n : \VY_n \to \VtildeY_n$ be the corresponding isomorphism.
In $\tilde Y_n$ the image $\alpha_n(Y_{n-1})$ is now replaced by $Y_n$
such that for any $x \in \VY_0$ the vertex $\psi_n \circ \dotsb \circ \psi_1(x)$
is identified with $\alpha_n(\psi_{n-1} \circ \dotsb \circ \psi_1 (x))$.
In the resulting graph we replace every $\theta$-complete graph $H$,
which is not completely contained in $Y_n$, by an isomorphic copy of $G$
such that any different vertices in the copy, that correspond to vertices in $B$,
are amalgamated with different vertices in $\VH$.
In general this procedure of replacing the $\theta$-complete graphs is not unique.
Let $Y_{n+1}$ be the resulting graph.

We set $F_{n+1} = F_n \cup \alpha_n(\VY_n \setminus \VY_{n-1})$. The function
\[ \psi_{n+1} : \VY_n \to F_{n+1} \]
defined by
\[ {\psi_{n+1}|}_{\VY_{n-1}} = \psi_n \qquad\text{and}\qquad
	{\psi_{n+1}|}_{\VY_n \setminus \VY_{n-1}} =
	{\alpha_n|}_{\VY_n \setminus \VY_{n-1}} \]
is an isomorphism between $Y_n$ and $(Y_{n+1})_{F_{n+1}}$.
The sequences $(Y_n)_{n\in\N}$ and $(F_n)_{n\in\N}$ are increasing
and $X = \bigcup_{n=1}^\infty Y_n$ is self-similar
with respect to $F = \bigcup_{n=1}^\infty F_n$ and $\psi : \VX \to F$
defined by $\psi|_{\VY_n} = \psi_{n+1}$.

This construction implies that (a) and (b) are satisfied. We have
\[ d_X(\psi(x),\psi(y)) \geq 2 d_X(x,y) \]
for the `blow-up' function $\psi$ and vertices $x$ and $y$ in $VX$ and
\[ d_X(\phi(v),\phi(w)) \leq 2 d_X(v,w)\]
for its inverse contraction $\phi=\psi^{-1}$ which is defined on $F$,
for details see \cite[Theorem~1]{kroen02green}.
This implies that a fixed-point $o$ has must lie in $\VK_\theta^i \cap B$
and this intersection is a singleton.
The set $\theta \{o\} = \VK_\theta^i \setminus \{o\}$ is connected in $X$,
thus $\VX \setminus \{o\}$ is connected.
\end{proof}

By choosing $\theta=2$ every connected finite graph
which is not complete can occur as cell graph of a self-similar graph satisfying (S1).
Any self-similar graph with a unique origin cell can be obtained
as a blow-up in the sense of Theorem~\ref{thm:blowup},
thus it characterises self-similar graphs
which remain connected by removing an origin vertex,
see Lemma~\ref{lem:banach}.

Let $X$ be a symmetrically self-similar graph and let $C$ be any 1-cell of $X$.
Then there exists exactly one 2-cell $C_2$, such that $\psi(C) \subset C_2$.
We denote the 1-cells being contained in $C_2$ by $C^1, \dotsc, C^\mu$.
The cell graph $\hat C$ consists of $\mu$ amalgamated $\theta$-complete graphs
$K^1_\theta, \dotsc, K^\mu_\theta$, where $\VK_\theta^i = \phi(\theta C^i)$
and $i \in \{1,\dotsc,\mu\}$, see \cite{kroen02growth}.
For a vertex $v \in \VhatC$ we call
\[ \beta(v) = \abs{\{i \in \{1,\dotsc,\mu\} \mid v \in \VK_\theta^i\}} \]
the \emph{branching number} of $v$.
By Axiom~(S2), $\Aut(\hat C)$ acts transitively on $\theta C$.
Thus all vertices $v$ of the boundary $\theta C$ have the same branching number $\beta(v)$.
Furthermore, Axiom~(S1) implies that all cell graphs are isomorphic
and $\beta(v)$ does also not depend on the choice of the cell $C$.
Consequently, we define the branching number $\beta$ of a symmetrically self-similar graph
as the branching number $\beta(v)$ of any vertex $v$ in $F$.
This branching number characterises symmetrically self-similar graphs with bounded geometry:

\begin{thm}\label{thm:bounded}
A symmetrically self-similar graph $X$ has bounded geometry if and only if $\beta=1$.
\end{thm}

This theorem is an immediate consequence of Theorem~3 in \cite{kroen02growth}.
The next theorem is a special case of Corollary~1 in \cite{kroen02growth}.

\begin{thm}\label{thm:unbounded}
Let $X$ be a symmetrically self-similar graph with unbounded geometry.
Then either there exists no origin vertex and $X$ is locally finite.
Or there is an origin vertex which is then the only vertex with infinite degree.
\end{thm}

In the following lemma we characterise those
self-similar graphs satisfying (S1), which are bipartite.

\begin{lem}\label{lem:bipartite}
Let $X$ be a self-similar graph
with respect to $F$ and $\psi$ satisfying Axiom~(S1).
Then $X$ is bipartite if and only if some cell graph is bipartite.
Furthermore, we have $\theta=2$ if $X$ is bipartite.
\end{lem}

\begin{proof}
If $X$ is bipartite then every cell graph is bipartite,
since it is a subgraph of $X$.

Now suppose that some cell graph is bipartite,
whence all cell graphs are bipartite thanks to Axiom~(S1).
A $\theta$-complete graph is bipartite if and only if $\theta \leq 2$.
Any cell graph $\hat C$ consists of $\mu$ amalgamated copies
of $\theta$-complete graphs, hence $\theta=2$.
Every closed path in $\hat C$ has even length.
Either all paths in $\hat C$ connecting the two vertices in $\theta C$ have
even length or they have all odd length.

Suppose that $X$ is not bipartite.
Then there exists a closed path $\pi$ in $X$ with minimal odd length.
This path cannot be completely contained in one cell graph,
hence it meets vertices in $F$.
Let $\pi = (x_0, \dots, x_n)$ with $x_0 = x_n \in F$
and $x_{k+n}=x_k$ for all $k \in \N$.
Let $x_i$, $x_j$ be vertices in $F$
such that $i<j$ and $x_k \notin F$ for all $i < k < j$.
Then $j < i+n$, since $\pi$ is not contained in one cell graph.
Moreover, we have $x_i \neq x_j$.
Otherwise $(x_i, \dotsc, x_j)$ would be a closed path in one cell graph
and $(x_j, \dotsc, x_{i+n})$ would be a closed path in $X$ with odd length
strictly smaller than $n$.
Let $0=i_0 < \dotsb < i_\ell=n$ be all those indices,
such that $x_{i_k} \in F$ for $k \in \{0,\dotsc,\ell\}$.
Since the lengths of the paths $(x_{i_{k-1}},\dotsc,x_{i_k})$
for $k \in \{1,\dotsc,\ell\}$ have same parity,
$(\phi(x_{i_0}),\dotsc,\phi(x_{i_\ell}))$ is a closed path in $X$
of odd length strictly smaller than $n$ in contradiction to the minimality of $\pi$.
\end{proof}

\section{Green function at the origin}\label{sec:green}

In this section we repeat and reformulate
results and definitions from \cite{kroen02green}. Let
\[ P = \bigl(p(x,y)\bigr)_{x,y \in \VX} \]
be the matrix of transition probabilities for the simple random walk
on a symmetrically self-similar graph $X$ which is locally finite.
The Green function for vertices $x$ and $y$ in $\VX$ is defined as the generating function
\[ G(x,y|z) = \sum_{n=0}^\infty p^{(n)}(x,y) z^n \]
of the $n$-step transition probabilities $p^{(n)}(x,y)$ from $x$ to $y$ for $z \in U(0,1)$,
where $U(z_0,r) = \{ z \in \C \mid \abs{z-z_0}<r \}$ for any $z_0 \in \C$ and $r > 0$.
If the Green function is considered as infinite dimensional matrix, we have
\[ \bigl(G(x,y|z)\bigr)_{x,y\in \VX} = \sum_{n=0}^\infty P^n z^n = (I-zP)^{-1}. \]

Let $C$ be a cell and let $B$ be a nonempty subset of $\theta C$. Then
\[ Q_B = \bigl(q_B(x,y)\bigr)_{x,y \in \VhatC} \]
denotes the transition matrix of the simple random walk
on $\hat C$ with absorbing boundary $B$. This means that
\[ q_B(x,y) = \frac{1}{\deg_{\hat C}(x)} \]
if $x \in \VhatC \setminus B$ is adjacent to a vertex $y$ in $\VhatC$
and $q(x,y)=0$ otherwise. We define
\[ Q_B^\ast(z) = \sum_{n=0}^\infty (zQ_B)^n = (I-zQ_B)^{-1}, \]
for all $z\in\C$ such that $I-zQ_B$ is invertible.
The entries of the matrix $Q_B^\ast(z)$ are rational functions in $z$
which are holomorphic on $U(0,r)$ for some $r>1$.
This follows from the fact that the simple random walk is absorbed by $B$ almost surely,
since the graph $\hat C$ is finite.

Let $v \in \theta C$, then we define the rational functions $d$ and $f$ by
\[ d(z) = \sum_{w \in \theta C \setminus \{v\}}
		\bigl(Q_{\theta C \setminus \{v\}}^\ast(z)\bigr)_{v,w}
	\qquad\text{and}\qquad
	f(z) = \bigl(Q_{\theta C \setminus \{v\}}^\ast(z)\bigr)_{v,v} \,. \]
The functions $d$ and $f$ are independent of the choice of $v$
because $\Aut(\hat C)$ acts transitively on $\theta C$.
Furthermore, double transitivity implies
\[ d(z) = (\theta - 1) \cdot \bigl(Q_{\theta C \setminus \{v\}}^\ast(z)\bigr)_{v,w} \]
for any vertex $w \in \theta C \setminus \{v\}$.
Now $d$ is the generating function of the probabilities
that the simple random walk on $\hat C$ starting in some vertex $v$ in $\theta C$
hits a vertex in $\theta C \setminus \{v\}$ for the first time after exactly $n$ steps.
Whereas $f$ is the generating function of the probabilities
that the random walk starting in $v$ returns to $v$ after $n$ steps
without hitting a vertex in $\theta C \setminus \{v\}$ before.
The start is counted as first visit, thus $f(0)=1$.
We call $d$ the \emph{transition function} and $f$ the \emph{return function} of $X$.

We write $G(z)$ instead of $G(o,o|z)$ if $X$ is a
self-similar graph with origin vertex $o$.

\begin{lem}\label{lem:funceq}
Let $X$ be a symmetrically self-similar graph with origin vertex $o$. Then
\begin{equation}\label{eq:funceq}
	G(z) = f(z) \cdot G(d(z))
\end{equation}
for all $z\in U(0,1)$.
\end{lem}

This functional equation was shown for the Sierpi\'nski graph in
\cite{rammal84random,rammal83random} and \cite{grabner97functional1}.
It is a special case of Lemma~6 in \cite{kroen02green}.

The basic idea behind this identity is the substitution $z \to d(z)$ which corresponds
to a combinatorial path substitution reflecting the self-similarity of the graph.
For more details see Section~4 in \cite{kroen02green}.
The number of orbits of $\Aut(\hat C)$ on the set $\{(v,w) \mid v,w \in \theta C\}$
is the number of `types' of transitions from one boundary vertex in $\theta C$ to another.
It is the number of variables that is needed
for the functional equations of the Green functions.
Similarly, this is the number of variables in the renormalization equations
for constructing Brownian motion on fractals (see for example \cite{lindstroem90brownian}).
Axiom~(S2) ensures that only one variable is needed.
Our techniques would apply to more general self-similar graphs
if the dynamics of multidimensional transition functions could be understood.

Let $\cd_0 \subset \C \cup \{\infty\}$ be the set of all poles of the function $f$,
and let $\cd$ be the set of all $z \in \C \cup \{\infty\}$
such that $d^n(z) \in \cd_0$ for some $n \in \N_0$.
The following lemma corresponds to the Lemmata~9, 11, 12 and Theorem~4 in \cite{kroen02green}.

\begin{lem}\label{lem:dynamics}
Let $X$ be a symmetrically self-similar graph with bounded geometry.
\begin{enum}
\item The point $z=1$ is a repelling fixed point of the transition function $d$.
\item The point $z=0$ is a super-attracting fixed point of $d$,
	and it is the only attracting fixed point.
\item The Fatou set $\cf$ of $d$ is the immediate basin of attraction of $z=0$.
\item The Julia set $\cj$ of $d$ is a subset of $(\R \cup \{\infty\}) \setminus (-1,1)$.
\item The set $\cd$ is a subset of $(\R \cup \{\infty\}) \setminus (-1,1)$.
\item Every Green function $G(x,y|\,\cdot\,)$ is holomorphic in $\cf \setminus \cd$
	and cannot be continued holomorphically to any point in $\cj$.
\end{enum}
\end{lem}

In \cite{kroen02green} this information was used for the description
of the spectrum of the Laplacian. Here it is the basic tool for studying
the transition probabilities.

We recall that a point $z_0$ is called \emph{exceptional point} of the map $d$
if the smallest completely invariant set containing $z_0$ is finite.
Let $z_0 \in \cf$ be an exceptional point of $d$.
Then the smallest completely invariant set containing $z_0$
must also contain $0$, since $\cf$ is the immediate basin of attraction of $0$.
Since a rational function of degree greater or equal two has at most two exceptional points,
we have $d(z_0)=0$ and $d(0)=z_0$,
see the proof of Theorem~4.1.2 in the book of Beardon \cite{beardon91iteration}.
Thus $d$ has at most the exceptional point $z=0$.
Since $0 \notin \cd_0$ the set of accumulation points of $\cd$ coincides with $\cj$,
see \cite[Theorem~4.2.7 and 4.2.8]{beardon91iteration}.
Hence $\overline{\cd} = \cd \cup \cj$ and
$\cf \setminus \cd = (\C \cup \{\infty\}) \setminus \overline{\cd}$.
We remark that either $z=0$ is exceptional and $d$ is conjugated to a polynomial
or there are no exceptional points at all. Actually, both cases can occur,
as we will show by examples in Sections~\ref{sec:line} and \ref{sec:examples}.

\begin{cor}\label{cor:continuation}
Let $X$ be a symmetrically self-similar graph with bounded geometry and origin vertex $o$.
The Green function $G$ has a unique holomorphic continuation to $\cf \setminus \cd$ and
\[ G(z) = \prod_{k=0}^\infty f(d^k(z)) \]
for all $z \in \cf \setminus \cd$.
This convergence is uniform on any compact subset of $\cf \setminus \cd$.
If $X$ is not bipartite then $\cf \setminus \cd$
contains an interval $(-r,-1]$ for some $r > 1$.
\end{cor}

\begin{proof}
The first part is a consequence of Lemma~\ref{lem:dynamics}.
Let $X$ be not bipartite. Then, by Lemma~\ref{lem:bipartite},
the cell graph $\hat C$ is not bipartite.
Hence $d$ is neither even nor odd. Therefore $\abs{d(-1)}<1$
and $\cf \setminus \cd$ contains $(-r,-1]$ for some $r>1$.
\end{proof}

\section{Geometric and probabilistic dimensions of self-similar graphs}\label{sec:geometric}

Let $X$ be a symmetrically self-similar graph with respect to $F$ and $\psi$.
We call \[ \tau = d'(1) \] the \emph{time scaling factor} of $X$.
Then $\tau$ is the expected number of steps of the simple random walk on the cell graph $\hat C$
for hitting $\theta C \setminus \{v\}$, when starting in some $v \in \theta C$.
Let $C_2$ be any 2-cell of $X$. Then the number of 1-cells in a 2-cell
\[ \mu = \abs{\{ C \in \cc(F) \mid C \subset \VhatC_2\}} \]
is called the \emph{mass scaling factor} of $X$.
Axiom (S1) implies that $\mu$ is independent of the choice of the 2-cell $C_2$.
Furthermore, we call
\[ d_s = \frac{2 \log(\mu)}{\log(\tau)} \]
the \emph{spectral dimension} of $X$. The spectral dimension
was first introduced in physics literature using the density of states of the Laplacian,
see \cite{alexander82density,rammal83random} and the references therein.
Later on, the exponent $d_s$ was studied from several points of view,
see \cite{barlow98diffusions,barlow88brownian,fukushima92spectral,
grabner97functional1,jones96transition,kigami01analysis,kigami93weyl}.

A measure $m$ on the set of vertices $\VY$ of a countable graph $Y$ is
\emph{invariant} with respect to a transition matrix
$P = (p(x,y))_{x,y \in \VY}$ if $m \cdot P = m$,
where $m$ is considered as row vector.
Written in coordinates this is
\[ m(x) = \sum_{y \in \VY} m(y) p(y,x) \]
for any $x \in \VY$. If $Y$ is connected and finite,
then the unique invariant probability measure of the simple random walk is given by
\[ m(x) = \frac{\deg_Y(x)}{\sum_{x \in \VY} \deg_Y(x)} = \frac{\deg_Y(x)}{2 \abs{\EY}}. \]

\begin{lem}\label{lem:invmeasure}
Let $X$ be a symmetrically self-similar graph,
and let $C$ be any cell of $X$. Then $m : \VhatC \to \R$ defined by
\[ m(x) = \frac{\beta(x)}{\mu \theta} \]
is the unique invariant probability measure of the simple random walk on $\hat C$.
\end{lem}

\begin{proof}
Let $x$ be any vertex in $\VhatC$.
Since the cell graph $\hat C$ consists of $\mu$ amalgamated $\theta$-complete graphs,
we have $\deg_{\hat C}(x) = (\theta-1) \beta(x)$ and
$2 \abs{\EhatC} = \mu \theta (\theta-1)$. This implies the lemma.
\end{proof}

In the following we will make use of path arguments.
For similar path decompositions see \cite{grabner97functional2},
\cite{grabner97functional1}, and especially \cite{kroen02green}.
These constructions can be seen in the context of more general generating functions,
see Goulden and Jackson \cite{goulden83combinatorial}.

Let $C$ be a cell of a symmetrically self-similar graph $X$.
For $v,w \in \theta C$ and $B \subset \theta C$
let $\Pi_B(v,w)$ be the set of paths from $v$ to $w$
which do not hit any vertex in $B$ except for their start and their end vertex.
The set of paths in $\Pi_B(v,w)$ with positive length
is denoted by $\Pi^\ast_B(v,w)$.

The \emph{weight} of a path $\pi = (x_0, \dotsc, x_n)$ in $\hat C$ is defined by
\[ W(\pi|z) = \prod_{i=0}^{n-1} \frac{z}{\deg_{\hat C}(x_i)} =
	\frac{z^n}{(\theta - 1)^n} \prod_{i=0}^{n-1} \frac{1}{\beta(x_i)} \]
for $z \in \C$. Now $W(\pi|1)$ is the probability,
that the simple random walk on $\hat C$, starting in $x_0$,
follows the path $\pi$ in its first $n$ steps.
For a set of paths $\Pi$ we set
\[ W(\Pi|z) = \sum_{\pi \in \Pi} W(\pi|z). \]
For $\Pi_1 \subset \Pi(x_1,x_2)$ and $\Pi_2 \subset \Pi(x_2,x_3)$
we write $\Pi_1 \circ \Pi_2$ for the set of all concatenations of
paths in $\Pi_1$ with paths in $\Pi_2$. Then
\[ W(\Pi_1 \circ \Pi_2|z) = W(\Pi_1|z) W(\Pi_2|z). \]

Let $v$ again be a vertex in $\theta C$.
We define $\hat u : U(0,1) \to \C$ and $u : U(0,1) \to \C$ by
\[ \hat u(z) = W(\Pi^\ast_{\{v\}}(v,v)|z) \qquad\text{and}\qquad
	u(z) = W(\Pi^\ast_{\theta C}(v,v)|z). \]
Then $\hat u$ is the generating function of the probabilities
that the simple random walk on $\hat C$,
starting in $v$, returns for the first time to $v$ after exactly $n$ steps, where $n>0$.
And $u$ is the generating function of the first return after $n$ steps
without hitting any vertex in $\theta C \setminus \{v\}$ before.
Since $\Aut(\hat C)$ acts transitively on $\theta C$,
the functions $\hat u$ and $u$ are independent of the choice of $v$.

\begin{lem}\label{lem:dfweight}
Let $C$ be a cell of a symmetrically self-similar graph $X$,
and let $v$ and $w$ be two different vertices in $\theta C$. Then
\[ f(z) = W(\Pi_{\theta C \setminus \{v\}}(v,v)|z) =
	1 +  W(\Pi^\ast_{\theta C \setminus \{v\}}(v,v)|z) =
	\frac{1}{1-u(z)}\]
and
\[ d(z) = (\theta - 1) \, W(\Pi_{\theta C \setminus \{v\}}(v,w)|z) =
	(\theta - 1) f(z) W(\Pi_{\theta C}(v,w)|z). \]
\end{lem}

\begin{proof}
Since $\Pi_{\theta C \setminus \{v\}}(v,v)$ is the set of all paths in $\hat C$ from $v$ to $v$
which do not hit vertices in $\theta C \setminus \{v\}$, we get
\[ f(z) = W(\Pi_{\theta C \setminus \{v\}}(v,v)|z) =
	1 +  W(\Pi^\ast_{\theta C \setminus \{v\}}(v,v)|z). \]
Every path in $\Pi^\ast_{\theta C \setminus \{v\}}(v,v)$ can be obtained by
concatenating a finite number of paths in $\Pi^\ast_{\theta C}(v,v)$, therefore
\[ f(z) = 1 + \sum_{k=1}^\infty u(z)^k = \frac{1}{1-u(z)}. \]
The definition of $d$ implies
\[ d(z) = \sum_{w \in \theta C \setminus \{v\}} W(\Pi_{\theta C \setminus \{v\}}(v,w)|z) =
	f(z) \sum_{w \in \theta C \setminus \{v\}} W(\Pi_{\theta C}(v,w)|z). \]
Here we used the fact, that
\[ \Pi_{\theta C \setminus \{v\}}(v,w) =
	\Pi_{\theta C \setminus \{v\}}(v,v) \circ \Pi_{\theta C}(v,w)). \]
Axiom (S2) now yields the rest of the statement.
\end{proof}

\begin{lem}\label{lem:firstreturn}
Let $X$ be a symmetrically self-similar graph, then
\[ \hat u(z) = u(z) + \frac{d(z)}{f(z)} \cdot \frac{1}{1-\frac{\theta-2}{\theta-1}d(z)}
			\cdot \frac{d(z)}{\theta-1}. \]
\end{lem}

\begin{proof}
Let $C$ be any cell of $X$ and let $v$ be any vertex in $\theta C$.
If we denote by $\Pi(v)$ those paths in $\Pi^\ast_{\{v\}}(v,v)$,
which hit a vertex in $\theta C \setminus \{v\}$, then
\[ \Pi^\ast_{\{v\}}(v,v) = \Pi^\ast_{\theta C}(v,v) \uplus \Pi(v) \]
and therefore $\hat u(z) = u(z) + W(\Pi(v)|z)$.
For an integer $n$ in $\N_0$ let $\Lambda_n(v)$ be the set of all $(n+1)$-tuples
$a = (a_0, \dots, a_n)$ in $(\theta C \setminus \{v\})^{n+1}$,
such that $a_{k} \neq a_{k+1}$ for $0\le k< n$.
If we write $\Pi_a(v)$ for the set
\[ \Pi_{\theta C \setminus \{a_0\}}(a_0,a_1) \circ \dotsb \circ 
	\Pi_{\theta C \setminus \{a_{n-1}\}}(a_{n-1}, a_n), \]
then, by Lemma~\ref{lem:dfweight}, 
\[ W(\Pi_a(v)|z) = \frac{d(z)^n}{(\theta-1)^n}. \]
Now we can decompose $\Pi(v)$ in the following way:
\[ \Pi(v) = \biguplus_{n \in \N_0} \biguplus_{a \in \Lambda_n(v)}
	\Pi_{\theta C}(v,a_0) \circ \Pi_a(v) \circ \Pi_{\theta C \setminus \{a_n\}}(a_n,v). \]
From this we obtain
\begin{align*}
W(\Pi(v)|z) &= \sum_{n \in \N_0} \sum_{a \in \Lambda_n(v)}
	W(\Pi_{\theta C}(v,a_0)|z) W(\Pi_a(v)|z) W(\Pi_{\theta C \setminus \{a_n\}}(a_n,v)|z) = \\
&= \sum_{n \in \N_0} \sum_{a \in \Lambda_n(v)} \frac{d(z)}{(\theta-1)f(x)} \cdot
	\frac{d(z)^n}{(\theta-1)^n} \cdot \frac{d(z)}{\theta-1} = \\
&= \frac{d(z)}{f(z)} \cdot \frac{1}{1-\frac{\theta-2}{\theta-1}d(z)}
			\cdot \frac{d(z)}{\theta-1}
\end{align*}
using Lemma~\ref{lem:dfweight} and the fact, that $\abs{\Lambda_n(v)} = (\theta-1) (\theta-2)^n$.
\end{proof}

For the following statements we recall that $\beta=1$
if and only if the graph has bounded geometry, see Theorem~\ref{thm:bounded}.
Let the simple random walk on $\hat C$ start in a vertex $v$ in $\theta C$.
Then $f(1)$ is the expected number of returns to $v$
(the start is counted as the first visit)
before hitting another vertex in $\theta C$.
The \emph{resistance scaling factor} is defined by the relation $\tau = \mu\rho$,
see \cite[Equation~(6.29)]{barlow98diffusions}.

\begin{thm}\label{thm:dim}
Let $X$ be a symmetrically self-similar graph. Then $f(1) = \beta\rho$.
\end{thm}

\begin{proof}
Since $\hat u(1)=1$ and $d(1)=1$, we obtain
\[ \hat u'(1) = \lim_{z \to 1} \frac{\hat u(z)-1}{z-1} =
	\lim_{z \to 1} \frac{d(z)+(\theta-1)}{f(z) ((\theta-1)-(\theta-2)d(z))}
		\cdot \frac{d(z)-1}{z-1} = \frac{\theta \tau}{f(1)} \]
using Lemma~\ref{lem:dfweight} and \ref{lem:firstreturn}.
Now $\hat u'(1)$ is the expected number of steps of the first return
of simple random walk on $\hat C$ starting in vertex $v$ in $\theta C$.
Thus $\hat u'(1) m(v) = 1$. A proof of this well known identity
can be found in the book of Br\'emaud, see \cite[Theorem~3.2]{bremaud99markov}.
By Lemma~\ref{lem:invmeasure}, we have
\[ \hat u'(1) = \frac{1}{m(v)} = \frac{\mu\theta}{\beta} \]
completing the proof.
\end{proof}

\begin{cor}
Let $X$ be a symmetrically self-similar graph. Then $\tau\beta > \mu$.
\end{cor}

If $X$ is a symmetrically self-similar graph with bounded geometry
then the last corollary implies $0 < d_s < 2$.

\section{Asymptotic analysis}\label{sec:asymptotic}

Throughout the rest of this paper let $X$ be a symmetrically self-similar graph
with bounded geometry and origin vertex $o$.
We recall that symmetrically self-similar graphs with an origin vertex
have either bounded geometry or
the origin vertex is the only vertex with infinite degree,
see Theorem~\ref{thm:unbounded}.

Since $d(1)=1$ and $\tau = d'(1) > 2$ there exists a holomorphic local inverse $d^{-1}$
of $d$ in a neighbourhood of $z=1$, which has a attracting fixed point at $z=1$.
For $n \in \N$ we write $d^{-n}$ for the $n$-fold iterate of $d^{-1}$.
Let $z_0 \in \C$, $r>0$ and $0 < \delta \leq \pi$; then we write $U_\delta(z_0,r)$
for the set of all $z \in \C \setminus \{z_0\}$,
such that $\abs{z-z_0}<r$ and $\abs{\arg(z_0-z)}<\delta$.
In particular, we have $U_\pi(1,r) = U(1,r) \setminus [1,1+r)$.
Furthermore, we write $S(b)$ for the horizontal strip
\[ S(b) = \bigl\{ z \in \C \bigm| \abs{\Im(z)} < \tfrac{\pi}{\log(b)} \bigr\}. \]

\begin{lem}\label{lem:period}
Let $H : U_\pi(1,R) \to \C$ be a holomorphic function for some $R>0$
which satisfies the equation $H(z) = H(d(z))$
whenever $z$ in $U_\pi(1,R)$ and $d(z) \in U_\pi(1,R)$.
Then there exists a $1$-periodic, holomorphic function $\omega$
defined on the strip $S(\tau)$ and a $r \in (0,R)$, such that
\[ H(z) = \omega\biggl(\frac{\log(1-z)}{\log(\tau)}\biggr) + \co_\delta(\abs{z-1}) \]
holds for $z \in U_\delta(1,r)$ for all $0 < \delta < \pi$.
The function $\omega$ is constant if and only if $H$ is constant.
\end{lem}

This lemma was proved by one of the authors in \cite[Lemma~5]{teufl01average}.
The proof of the lemma and the following theorem are based on ideas of
de Bruijn \cite{bruijn79asymptotic}, Odlyzko \cite{odlyzko82periodic} and
Grabner and Woess \cite{grabner97functional1}. A detailed discussion of the methods
can be found in the book \cite[Section~16]{woess00random}.

The point $z=1$ is a non-polar singularity of any Green function
of any recurrent random walk on an infinite, locally finite graph,
see for example \cite[Lemma~10]{kroen02green}. Lemma~3 in \cite{kroen02green} implies
that the simple random walk on any symmetrically self-similar graph
with bounded geometry is recurrent.

\begin{thm}\label{thm:greenasymp}
Let $X$ be a symmetrically self-similar graph with bounded geometry and origin vertex $o$,
then there exists a $1$-periodic, holomorphic function $\omega$ on the strip $S(\tau)$,
such that the Green function $G$ has the local singular expansion
\[ G(z) = (1-z)^{\frac{d_s}{2}-1}
	\biggl(\omega\biggl(\frac{\log(1-z)}{\log(\tau)}\biggr) + \co_\delta(\abs{z-1})\biggr) \]
for $z \in U_\delta(1,r)$ and $0 < \delta < \pi$, where $r > 0$.
\end{thm}

\begin{proof}
We substitute $G(z) = (1-z)^\alpha H^+(z)$ into \eqref{eq:funceq},
where $H^+$ is holomorphic in $U_\pi(1,1)$
and $\alpha \in \R$ is some constant. This yields
$H^+(z) = f(z) q(z)^\alpha H^+(d(z))$, where $q$ is given by
\[ q(z) = \frac{1-d(z)}{1-z} \]
and $q(1) = \tau$. Furthermore we define
\[ H^-(z) = \prod_{n=-\infty}^{-1} f(d^n(z)) q(d^n(z))^\alpha. \]
Since $f$ and $q$ are holomorphic in $U(1,R)$ for some $R > 0$,
the last product converges if $f(1) q(1)^\alpha = 1$.
Hence, we have to choose
\[ \alpha = \frac{\log(\mu)}{\log(\tau)} - 1 = \frac{d_s}{2} - 1. \]
Moreover, the convergence is uniform in $U(1,R)$,
and therefore the function $H^-$ is holomorphic in $U(1,R)$.
Obviously, $H^-(d(z)) = H^-(z) f(z) q(z)^\alpha$ holds
for all $z \in U(1,R)$ such that $d(z) \in U(1,R)$.
Let $H : U_\pi(1,R) \to \C$ be defined by $H(z) = H^-(z) H^+(z)$,
then $H$ is holomorphic and satisfies
\begin{align*}
H(z) &= H^-(z) H^+(z) = H^-(z) f(z) q(z)^\alpha H^+(d(z)) \\
	&= H^-(d(z)) H^+(d(z)) = H(d(z)),
\end{align*}
whenever $z$ and $d(z)$ are contained in $U_\pi(1,R)$.
Applying Lemma~\ref{lem:period} to $H$ we get a $1$-periodic holomorphic function $\omega$
on the strip $S(\tau)$ and a $r \in (0,R)$, such that
\[ H(z) = \omega\biggl(\frac{\log(1-z)}{\log(\tau)}\biggr) + \co_\delta(\abs{z-1}) \]
holds for $z \in U_\delta(1,r)$ for all $0 < \delta < \pi$.
Since $H^-(1) = 1$ and $H^-$ is holomorphic in $U(1,r)$, we finally obtain
\[ G(z) = (1-z)^\alpha \frac{H(z)}{H^-(z)} = (1-z)^\alpha
	\biggl(\omega\biggl(\frac{\log(1-z)}{\log(\tau)}\biggr) + \co_\delta(\abs{z-1})\biggr) \]
for $z \in U_\delta(1,r)$ and $0 < \delta < \pi$.
\end{proof}

We will now use the method of `singularity analysis' due to
Flajolet and Odlyzko \cite{flajolet90singularity},
in order to obtain the asymptotic behaviour of $p^{(n)}(o,o)$.
It should be mentioned here that a real Tauberian
theorem would not reveal the oscillating nature.

\begin{thm}\label{thm:asymp}
Let $X$ be a symmetrically self-similar graph with bounded geometry and origin vertex $o$.
If $X$ is not bipartite, then
\[ p^{(n)}(o,o) = n^{-\frac{d_s}{2}}
	\biggl( \sigma\biggl(\frac{\log(n)}{\log(\tau)}\biggr) + \co(n^{-1})\biggr), \]
where $\sigma : S(\tau^2) \to \C$ is a $1$-periodic, holomorphic function
given by its Fourier series
\[ \sigma(z) = \sum_{k=-\infty}^\infty
	\frac{\hat\omega(-k)}{\Gamma\bigl(1-\frac{\log(\mu)+2k\pi i}{\log(\tau)}\bigr)} \cdot
	e^{2k\pi i z}. \]
If $X$ is bipartite, then $p^{(2n+1)}(o,o)=0$ and
\[ p^{(2n)}(o,o) = 2(2n)^{-\frac{d_s}{2}}
	\biggl( \sigma\biggl(\frac{\log(2n)}{\log(\tau)}\biggr) + \co(n^{-1})\biggr). \]
\end{thm}

\begin{proof}
Since the $1$-periodic function $\omega$ of Theorem~\ref{thm:greenasymp}
is holomorphic in $S(\tau)$, we have
\[ \hat\omega(k) = \co_\eps\Bigl(e^{-2\pi \abs{k} (\frac{\pi}{\log(\tau)}-\eps)}\Bigr) \]
for all $k \in \Z$ and $\eps > 0$,
where $\hat\omega(k)$ denotes the $k$-th Fourier coefficient of $\omega$.
Therefore it is possible to apply the method of `singularity analysis'.
For the technical details of the transfer we refer to \cite[Section~5]{odlyzko82periodic}
and \cite[Section~16]{woess00random}.

If $X$ is not bipartite, then $z=1$ is the only singularity of $G$ at the boundary of the unit disk.
The asymptotic of $p^{(n)}(o,o)$ is now a consequence of the singularity analysis.

If $X$ is bipartite, then $G$ is an even function. Hence $p^{(2n+1)}(o,o)=0$.
Furthermore, $G$ has the same singularity at $z=-1$ as at $z=1$.
Thus we have to add up the terms corresponding to the singularities $z=1$ and $z=-1$
for the asymptotics of $p^{(2n)}(o,o)$.
\end{proof}

We remark here that the function $H$ of the proof of Theorem~\ref{thm:greenasymp}
maps the interval $(1-r,1)$ to $\R$. By inspection
of the proof of Lemma~\ref{lem:period} in \cite{teufl01average},
we see that both $\omega$ and $\sigma$ map $\R$ to $\R$.
Hence the Fourier coefficients of $\omega$ and $\sigma$ satisfy
\[ \hat\omega(-k) = \overline{\hat\omega(k)}
	\qquad\text{and}\qquad
	\hat\sigma(-k) = \overline{\hat\sigma(k)} \]
for all $k \in \Z$. It is possible to compute the
Fourier coefficients of $\omega$ for a given graph $X$ numerically,
see \cite{teufl01average}.

\section{Oscillation}\label{sec:oscillation}

Let $X$ be a symmetrically self-similar graph with bounded geometry and origin vertex $o$.
We have shown that the asymptotics of $p^{(n)}(o,o)$ carries
an oscillating factor given by
\[ \sigma\biggl(\frac{\log(n)}{\log(\tau)}\biggr), \]
where $\sigma$ is a $1$-periodic, holomorphic function given
by the Fourier series in Theorem~\ref{thm:asymp}.
In this section we will give a necessary condition for the case
that the $1$-periodic function $\sigma$ is constant.
First of all we note that the following statements are equivalent:
\begin{itemize}
\item The $1$-periodic function $\sigma$ is constant.
\item The $1$-periodic function $\omega$ of Theorem~\ref{thm:greenasymp} is constant.
\item The function $H$ given in the proof of Theorem~\ref{thm:greenasymp} is constant.
\end{itemize}

A closed interval in $\R \cup \{\infty\}$ is denoted by $[a,b]$,
where $a, b \in \R \cup \{\infty\}$. Here $a$ and $b$ do not
necessarily satisfy the inequality $a \leq b$. For example, $[1,-2]$ is the set
$\{z \in \R \mid z \geq 1\} \cup \{\infty\} \cup \{z \in \R \mid z \leq -2\}$.
The same notation is used for open, respectively half open, intervals.

\begin{thm}\label{thm:oscillation}
If $\sigma$ is constant then the Julia set $\cj$ of the transition function $d$
is a closed interval $[1,a]$ in the extended real line $\R \cup \{\infty\}$, where $a \in (1,-1]$.
\end{thm}

\begin{proof}
The Julia set $\cj \subset \R \cup \{\infty\}$ is either a Cantor set
or a closed interval, see \cite[Theorem~5.7.1]{beardon91iteration}.
We will use the notation of the proof of Theorem~\ref{thm:greenasymp}.
If $\sigma$ is constant then the function $H$ is also constant.
Hence we obtain
\[ G(z) = \frac{C (1-z)^\alpha}{H^-(z)} \]
for $z \in U(1,r)$, where $C$ is some constant and $r > 0$.
Since $z=1$ is not a polar singularity by Lemma~10 in \cite{kroen02green},
we have $\alpha \notin \Z$.
Therefore $G$ has no analytic continuation to the interval $[1,1+r)$
and Corollary~\ref{cor:continuation} implies that $[1,1+r) \subset \cj \cup \cd$.
Furthermore, the set of accumulation points of $\cd$ coincide with $\cj$.
Thus $\cj$ contains the interval $[1,1+r)$
and must therefore be itself a closed interval in $\R \cup \{\infty\}$.
\end{proof}

\begin{cor}\label{cor:oscillation}
If $\cj$ is a Cantor set, then $\sigma$ is non-constant.
\end{cor}

It is not clear, if the opposite direction of this corollary is also true.
For rational functions with real Julia set there is a general method to decide
whether their Julia set is an interval or a Cantor set.
Therefore we refer to Inninger's detailed discussion of rational functions
in Section~3 and Section~4 of \cite{inninger01rational}.
However, for a given transition function $d$
the following observation may be useful as well:
If we can find a real number $z_0 > 1$
which is smaller than some $d$-backwards iterate of $z=1$
together with an integer $n$ such that $\abs{d^n(z_0)}<1$,
then $\cj$ set is a Cantor set.

\section{Transition probabilities between arbitrary vertices}\label{sec:transition}

Let $X$ be a symmetrically self-similar graph with bounded geometry
and let $P$ be the matrix of transition probabilities for the simple random walk on $X$.

\begin{lem}
If $X$ is not bipartite, then $P$ is strongly aperiodic, that is,
there exists a number $n_0$ such that
$\inf\{p^{(n)}(x,x) \mid x \in X\} > 0$ for all $n \geq n_0$.
In general, $P^2$ is strongly aperiodic.
\end{lem}

\begin{proof}
Since $X$ has bounded geometry and is connected, there exists a $c_2 > 0$
such that $p^{(2)}(x,x) \geq c_2$ for all $x \in \VX$.
Thus $P^2$ is strongly aperiodic.

If $X$ is not bipartite, then the cell graph $\hat C$ of $X$
is not bipartite by Lemma~\ref{lem:bipartite}.
Hence for any vertex $v \in \VhatC$ there exists
a closed path in $\hat C$ from $v$ to $v$ of odd length.
By adding a path of the form $(v,w,v)$,
where $w$ is adjacent to $v$ in $\hat C$,
we can extend such a path for $v$ by an even number of edges.
Thus there is an odd number $\ell$ such that for any
$x \in \VhatC$ we can find a closed path from $v$ to $v$ of length $\ell$.
Since any vertex $x \in \VX$ is contained in at least one cell graph,
there is a real number $c_1 > 0$ such that $p^{(\ell)}(x,x) \geq c_1$ for all $x \in \VX$.
As every $n \geq \ell$ can be written in the form $r \cdot 2 + s \cdot \ell$ with $r,s \in \N_0$,
strong aperiodicity follows.

We remark here, that it is possible to choose $\ell=3$, if $\theta \geq 3$,
since the cell graph consists of $\mu$ copies of the $\theta$-complete graph.
\end{proof}

\begin{thm}\label{thm:offdiagonal}
Let $x$ and $y$ be two vertices of a symmetrically self-similar graph $X$
with bounded geometry and origin vertex $o$.
If $X$ is not bipartite, then 
\[ p^{(n)}(x,y) = n^{-\frac{d_s}{2}}
	\biggl( \sigma\biggl(\frac{\log(n)}{\log(\tau)}\biggr) + o(1)\biggr). \]
If $X$ is bipartite, let $r$ be an element of $\{0,1\}$ such that $d(x,y) \equiv r \bmod{2}$.
Then $p^{(2n+1-r)}(x,y)=0$ and
\[ p^{(2n+r)}(x,y) = 2(2n+r)^{-\frac{d_s}{2}}
	\biggl( \sigma\biggl(\frac{\log(2n+r)}{\log(\tau)}\biggr) + o(1)\biggr). \]
\end{thm}

\begin{proof}
If $X$ is not bipartite then
\[ \lim_{n \to \infty} \frac{p^{(n)}(x,y)}{p^{(n)}(y,y)} = 1, \]
by \cite[Theorem~2]{grabner97functional1}, since $P$ is
irreducible, strongly aperiodic and recurrent. Thus
\[ \frac{p^{(n)}(x,y)}{p^{(n)}(o,o)} = \frac{p^{(n)}(x,y)}{p^{(n)}(y,y)} \cdot
	\frac{p^{(n)}(y,y)}{p^{(n)}(y,o)} \cdot \frac{p^{(n)}(y,o)}{p^{(n)}(o,o)} = 1 + o(1). \]
Now Theorem~\ref{thm:asymp} implies the statement.

If $X$ is bipartite then we can partition $\VX$ into classes $V_1$ and $V_2$,
such that edges in $\EX$ only connect vertices in $V_1$ with vertices in $V_2$.
Then $P^2$ restricted to $V_i$ for $i \in \{1,2\}$
is irreducible, strongly aperiodic and recurrent.
If $r=0$, then $x$ and $y$ are in the same class.
Hence $p^{(2n+1)}(x,y)=0$ and
\[ \frac{p^{(2n)}(x,y)}{p^{(2n)}(o,o)} = 1 + o(1). \]
If $r=1$, then $p^{(2n)}(x,y)=0$ and
\[ \frac{p^{(2n+1)}(x,y)}{p^{(2n)}(o,o)} =
	\sum_{v \in \VX} \frac{p(x,v) p^{(2n)}(v,x)}{p^{(2n)}(o,o)} =
	\sum_{v \in \VX} p(x,v)(1 + o(1)) = 1 + o(1). \]
Note that the above sum is finite, since there are
only finitely many vertices $v \in \VX$ such that $p(x,v) > 0$. Furthermore,
\[ (2n)^{-\frac{d_s}{2}} = (2n+1)^{-\frac{d_s}{2}} (1 + o(1)). \]
As $\sigma$ is a $1$-periodic, holomorphic function, it is uniformly continuous on the real line.
We have $\log(2n+1) - \log(2n) = O(n^{-1})$ and therefore
\[ \sigma\biggl(\frac{\log(2n)}{\log(\tau)}\biggr) =
	\sigma\biggl(\frac{\log(2n+1)}{\log(\tau)}\biggr) + o(1) \]
which completes the proof.
\end{proof}

Jones computed estimates for the transition probabilities
$p^{(n)}(x,y)$ on the Sierpi\'nski graph
which are uniform in space ($x,y \in \VX$) and time ($n \in \N$).
His result yields the correct asymptotic type
\[ p^{(n)}(x,y) \asymp n^{-\frac{\log(3)}{\log(5)}} \]
as $n \to \infty$, see \cite{jones96transition}.

\section{Geodesic transition probabilities and the line graph}\label{sec:line}

\begin{exmp}\label{exmp:line}
Let the cell graphs of a self-similar graph $X$ be isomorphic
to the line $L_n$ of length $n \geq 2$ (Figure~\ref{fig:line}a).
If the graph $X$ has an origin vertex $o$, then it consists of an arbitrary number
of one-sided infinite lines which are amalgamated in their initial vertex $o$,
see Figure~\ref{fig:line}b. If there is no origin vertex then $X$ is the two-sided infinite line.
\begin{figure}[ht]
\begin{picture}(290,80)(0,-35)
\refstepcounter{fig}\label{fig:line}
\put(0,5){%
\multiput(0,0)(30,0){3}{\circle*{2.5}}
\dashline{2}(60,0)(80,0)
\dashline{2}(110,0)(130,0)
\drawline(130,0)(160,0)
\multiput(130,0)(30,0){2}{\circle*{2.5}}
\drawline(0,0)(60,0)
\put(-3,-10){$x_0$}
\put(27,-10){$x_1$}
\put(57,-10){$x_2$}
\put(127,-10){$x_{n-1}$}
\put(157,-10){$x_n$}
\put(80,10){$L_n$}}
\put(240,14){%
\setlength{\unitlength}{0.2pt}
\put(100,0){\circle*{11}}
\put(200,0){\circle*{11}}
\put(93.97,34.2){\circle*{11}}
\put(187.94,68.4){\circle*{11}}
\put(76.6,64.28){\circle*{11}}
\put(153.2,128.56){\circle*{11}}
\put(93.97,-34.2){\circle*{11}}
\put(187.94,-68.4){\circle*{11}}
\put(0,0){\circle*{11}}
\put(-40,-10){$o$}
\put(30,100){$X$}
\drawline(187.94,68.4)(0,0)(153.2,128.56)(0,0)(200,0)(0,0)(187.94,-68.4)
\dashline[20]{10}(187.94,68.4)(244.322,88.92)
\dashline[20]{10}(153.2,128.56)(199.16,167.128)
\dashline[20]{10}(187.94,-68.4)(244.322,-88.92)
\dashline[20]{10}(0,0)(199.16,-167.128)
\dashline[20]{10}(200,0)(260,0)}
\put(60,-35){\emph{Figure}~\ref{fig:line}a}
\put(238,-35){\emph{Figure}~\ref{fig:line}b}
\end{picture}
\end{figure}
The corresponding transition function $d_n$ is conjugated to the $n$-th Chebychev polynomial
with conjugacy map $z\to \smash{\frac{1}{z}}$,
\[ d_n(z) = z^n \, \Biggl( \sum_{k=0}^{\lfloor \frac{n}{2}\rfloor}
				\binom{n}{2k} (1-z^2)^k \Biggr)^{-1}
	\qquad\text{and}\qquad d_n(z) = T_n(\tfrac{1}{z})^{-1}, \]
see \cite[Example~1]{kroen02green}.
As $\cj(T_n)=[-1,1]$ we get $\cj(d_n) = [1,-1] = (\R \cup \{\infty\}) \setminus (-1,1)$.
We have $\theta=2$, $\mu=n$, $\beta=1$, $\rho=n$ and $\tau=n^2$.
Furthermore, if there is an origin vertex, we have
\[ G(z) = \frac{1}{\sqrt{1-z^2}} \]
for $z \in \C \setminus \cj(d_n)$. Hence the function $\sigma$ is constant.
\end{exmp}

To prove that lines as cell graphs are characterised by its return and
transition functions we state the following lemma.

\begin{lem}\label{lem:lines}
Let $Y$ be a locally finite, connected graph,
and let $q^{(n)}(x,y)$ be the $n$-step transition probabilities of
a reversible Markov chain on $Y$ which is of nearest neighbourhood type.
Then $q^{(n)}(x,y) q^{(n)}(y,x) \leq 4^{1-n}$
for different vertices $x$ and $y$, which are at distance $n=d(x,y)$ from each other.
Here equality holds, if and only if
\[ \sum_{\substack{w \in \VY\\d(y,w)<d(y,x)}} \mkern -20mu q(x,w) \; = \mkern -10mu
	\sum_{\substack{w \in \VY\\d(x,w)<d(x,y)}} \mkern -20mu q(y,w) \; = \; 1 \]
and
\[ \sum_{\substack{w \in \VY\\d(x,w)<d(x,v)}} \mkern -20mu q(v,w) \; = \mkern -10mu
	\sum_{\substack{w \in \VY\\d(y,w)<d(y,v)}} \mkern -20mu q(v,w) \; = \; \frac{1}{2} \]
for all $v \in \VY \setminus \{x,y\}$,
this is, the Markov chain can be projected to
the simple random walk on the line $L_n$.
\end{lem}

\begin{proof}
We prove a more general statement which is more adequate for induction:
For a vertex $v \in \VY$ and a finite set of vertices $W$ of $\VY$ we write
$d(v,W) = \min\{d(v,w) \mid w \in W\}$.
Let $n \in \N$ and let $A$ and $B$ be finite subsets of $\VY$
such that $d(a,B)=d(b,A)=n$ for all $a \in A$ and $b \in B$.
For $k \in \{0,\dotsc,n\}$ we introduce the level set $D_k \subset \VY$ by
\[ D_k = \{ x \in \VY \mid d(x,A)=k, \, d(x,B)=n-k \}. \]
These level sets are finite, $D_0=A$ and $D_n=B$.
Furthermore, we define $D_{-1}=D_{n+1}=\emptyset$, and
\[ m_-(x) = m(x) \sum_{y \in D_{k-1}} q(x,y), \qquad
	m_+(x) = m(x) \sum_{y \in D_{k+1}} q(x,y), \]
$m_0(x) = m(x) - m_-(x) - m_+(x)$ for a vertex $x$ in $D_k$,
where $k \in \{0, \dotsc, n\}$. We define
\[ Q(D_i,D_j) = \sum_{v \in D_i}\sum_{w \in D_j}
	\frac{m(v)q^{(n)}(v,w)}{m_+(v)} \cdot
	\frac{m(w)q^{(n)}(w,v)}{m_-(w)}. \]
for $i,j \in \{0,\dotsc,n\}$. Then we claim that $Q(A,B) \leq 4^{1-n}$.
Equality holds, if and only $m_-(x) = m_+(x) = \frac{1}{2} m(x)$
for any vertex $x$ in $D_k$, where $k \in \{1,\dotsc,n-1\}$.

First of all we notice that the condition for equality
indeed imply $Q(A,B)=4^{1-n}$.

We use induction over $n$. For $n=1$ the inequality is immediate.
Now let $n>1$ and let the statement hold for all $k<n$.
Factorizing with respect to the $k$-th step and using the reversibility yields
\[ Q(A,B) = \sum_{a\in A}\sum_{b\in B} \frac{1}{m_+(a)m_-(b)} \,
	\biggl( \sum_{x\in D_k} m(x) q^{(k)}(x,a) q^{(n-k)}(x,b) \biggr)^2 \]
for $k \in \{1,\dotsc,n-1\}$. Now we use the inequality
$4m_-(x)m_+(x) \leq m(x)^2$ for $x \in D_k$, where $k \in \{1,\dotsc,n-1\}$.
Here equality holds if and only if $m_-(x) = m_+(x) = \frac{1}{2} m(x)$.
Thus we obtain
\begin{align*}
Q(A,B) &\leq \frac{1}{4} \sum_{a\in A}\sum_{b\in B} \frac{1}{m_+(a)m_-(b)} \,
	\biggl( \sum_{x\in D_k} \frac{m(x)q^{(k)}(x,a)}{\sqrt{m_-(x)}} \cdot
		\frac{m(x)q^{(n-k)}(x,b)}{\sqrt{m_+(x)}} \biggr)^2 \\
&\leq \frac{1}{4} \sum_{a\in A}\sum_{b\in B} \frac{1}{m_+(a)m_-(b)} \,
	\biggl( \sum_{x\in D_k} \frac{m(x)^2 q^{(k)}(x,a)^2}{m_-(x)} \biggr)
	\biggl( \sum_{x\in D_k} \frac{m(x)^2 q^{(n-k)}(x,b)^2}{m_+(x)} \biggr)
\end{align*}
using Cauchy-Schwarz inequality for the second step.
Now reversibility together with the induction hypothesis
for $Q(A,D_k)$ and $Q(D_k,B)$ implies that
\[ Q(A,B) \leq \tfrac{1}{4}\,Q(A,D_k)\,Q(D_k,B) \leq 4^{1-n}. \]
As we can factorize with respect to any $k \in \{1,\dotsc,n-1\}$ we get
the postulated conditions for equality.
\end{proof}

\begin{cor}\label{cor:lines}
Let $d_n$ be the transition function and let $f_n$ be the return function
of a cell graph which is isomorphic to the line $L_n$ of length $n \geq 2$.
Then there is no other cell graph of a self-similar graph with bounded geometry
with transition function $d_n$ or return function $f_n$.
\end{cor}

\begin{proof}
Let the transition function of a cell graph be the function $d_n$.
Then $d_n$ is either even or odd, whence the cell graph must be bipartite.
By Lemma~\ref{lem:bipartite} the self-similar graph itself must be bipartite
and therefore $\theta=2$. As we have bounded geometry, $\beta=1$ by Theorem~\ref{thm:bounded}.
The series expansion of $d_n$ in $z=0$ starts with $2^{1-n} z^n$.
Now Lemma~\ref{lem:lines} implies that the cell graph $\hat C$ must be the line $L_n$.

Suppose $f_n$ is the return function of a cell graph $\hat C$
which is not isomorphic to $L_n$ and let $v$ be a vertex in the boundary $\theta C$.
Since $f_n$ is even, we have $\theta=2$ and $\beta=1$ as before.
There are vertices in $C$ which have edge degree strictly larger than 2 in the graph $\hat C$.
Let $m$ be the minimal distance of $v$ to a vertex $x \in C$ with $\deg_{\hat C}(x) \geq 3$.
Then the probability of returning to $v$ in exactly $2m$ steps is
strictly smaller than the corresponding probability for the line $L_n$.
This means that the $2m$-th coefficients of the series expansions
of the return functions of $\hat C$ and $L_n$ around $z=0$ are different.
\end{proof}

\section{Further Examples}\label{sec:examples}

Well known Examples of symmetrically self-similar graphs are
the Vi\v cek graph and the Sierpi\'nski graph.
Both graphs do not belong to the following simple example
of an infinite class of self-similar graphs.

\begin{exmp}\label{exmp:flake}
Let the cell graph $\hat{C}^n$ consist of $n$ copies of the $n$-complete graph
which are amalgamated in one vertex $x$, where $n \geq 2$.
The boundary of the cell contains exactly one vertex from each of these copies,
but not the vertex $x$. Starting with these cell graphs
we can construct different symmetrically self-similar graphs
in the sense of Theorem~\ref{thm:blowup}.
The four vertices of the cell graph $\hat{C}^4$ in Figure~\ref{fig:flake4}a
which constitute the boundary $\theta C^4$ of the cell are drawn fat.
Figure~\ref{fig:flake4}b shows the corresponding 4-cell graph $\hat{C}^4_4$.
\begin{figure}[ht]
\begin{picture}(345,235)
\refstepcounter{fig}\label{fig:flake4}
\put(0,85){%
\put(0,0){\circle*{4}}
\put(40,40){\circle*{2}}
\put(80,80){\circle*{4}}
\put(80,0){\circle*{4}}
\put(0,80){\circle*{4}}
\put(5,45){\circle*{2}}
\put(5,35){\circle*{2}}
\put(75,35){\circle*{2}}
\put(75,45){\circle*{2}}
\put(35,5){\circle*{2}}
\put(35,75){\circle*{2}}
\put(45,5){\circle*{2}}
\put(45,75){\circle*{2}}
\put(0,0){\line(1,1){80}}
\put(0,80){\line(1,-1){80}}
\drawline(5,35)(75,45)(80,80)(45,75)(35,5)(0,0)(5,35)(35,5)
\drawline(75,35)(5,45)(0,80)(35,75)(5,45)
\drawline(35,75)(45,5)(80,0)(75,35)(45,5)
\drawline(45,75)(75,45)
\put(20,-20){\emph{Figure}~\ref{fig:flake4}a}}
\put(120,0){%
\put(0,10){\includegraphics*[width=8cm]{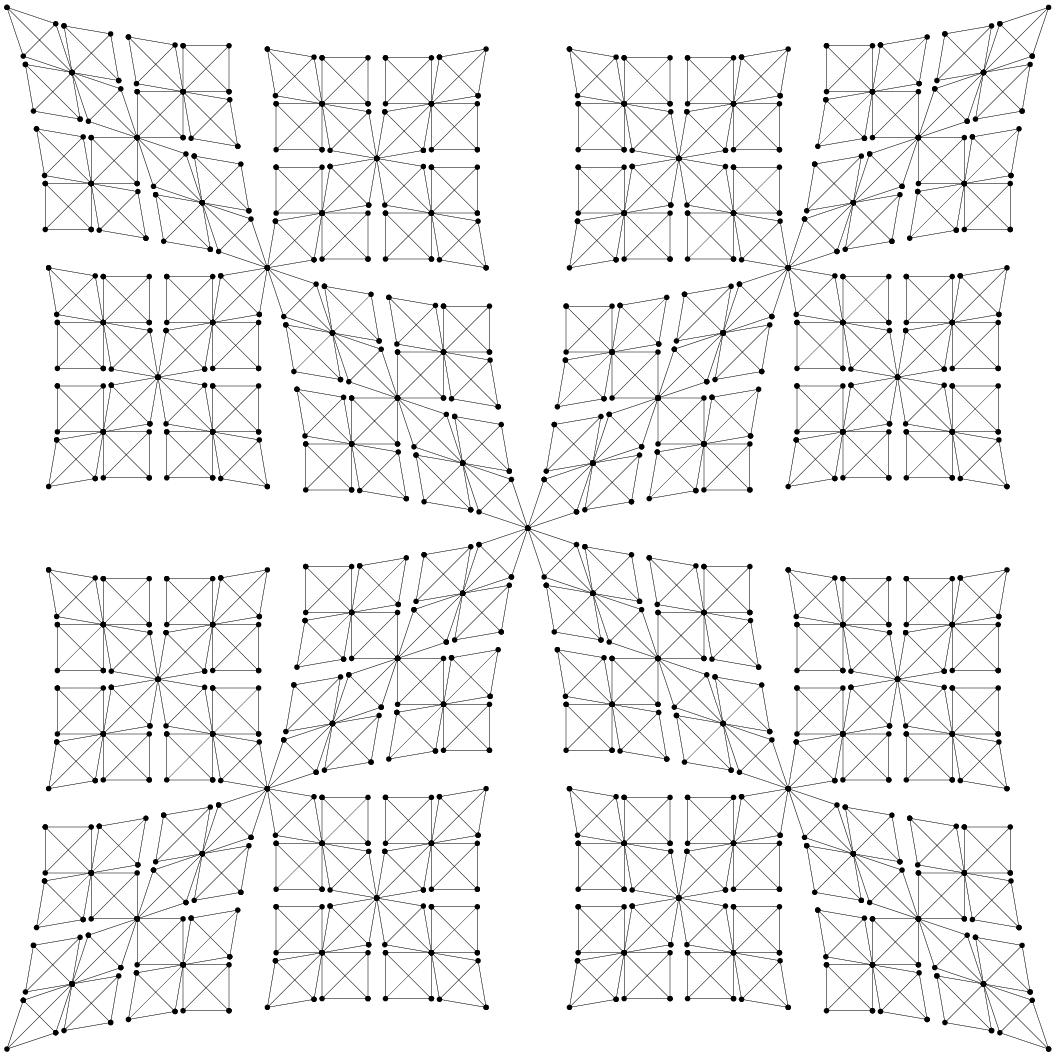}}
\put(90,0){\emph{Figure}~\ref{fig:flake4}b}}
\end{picture}
\end{figure}
The transition function $d_n$ of $\hat{C}^n$ is given by
\[ d_n(z) = \frac{z^2}{n(n-1) - 2n(n-2)z + (n^2-3n+1)z^2} \]
for $n \geq 2$. For the parameters we have $\mu=\theta=n$, $\beta=1$, $\tau=2n$ and $\rho=2$.
We remark that for $n=2$ the cells are isomorphic to the line $L_2$,
$d_2(z) = T_2(\frac{1}{z})^{-1}$, and $\cj(d_2) = [1,-1]$, see Example~\ref{exmp:line}.
For $n \geq 3$ we have
\[ d\biggl(\frac{n-1}{n-3}\biggr) = 1, \qquad\text{but}\qquad
	d^2\biggl(\frac{n-1}{n-2}\biggr) = \frac{n-1}{n(n-1)-1} \in (0,1). \]
Here $\smash{\frac{n-1}{n-3}}$ has to be read as $\infty$ if $n=3$.
As $1$ and $\smash{\frac{n-1}{n-3}}$ belongs to $\cj(d_n)$,
whereas $0$ and $\smash{\frac{n-1}{n-2}}$ does not,
the Julia set is not connected in $\R \cup \{\infty\}$.
Hence $\cj(d_n)$ is a Cantor set for $n \geq 3$.
In the case that there is an origin vertex
we can conclude that the $1$-periodic function $\sigma$ is not constant,
using Corollary~\ref{cor:oscillation}.
\end{exmp}

Many transition functions are of the form $d(z)=z^n/P(z)$,
where $n$ is the minimal number of steps from one point
in the boundary of a cell $C$ to another (in other words $n=\diam(\theta C)$)
and $P$ is a polynomial with integer coefficients.
As the following example shows, this is not true for all transition functions.

\begin{exmp}\label{exmp:34sierp}
The 3-dimensional, 4-scaled Sierpi\'nski graph consists of 20 amalgamated 4-complete graphs.
The 1-cell and 2-cell can be seen in Figure~\ref{fig:34sierp}.
Their boundary vertices are drawn fat.
Each 4-complete graph is represented by a solid tetrahedron.
\begin{figure}
\refstepcounter{fig}\label{fig:34sierp}
\begin{center}
\includegraphics*[width=5cm]{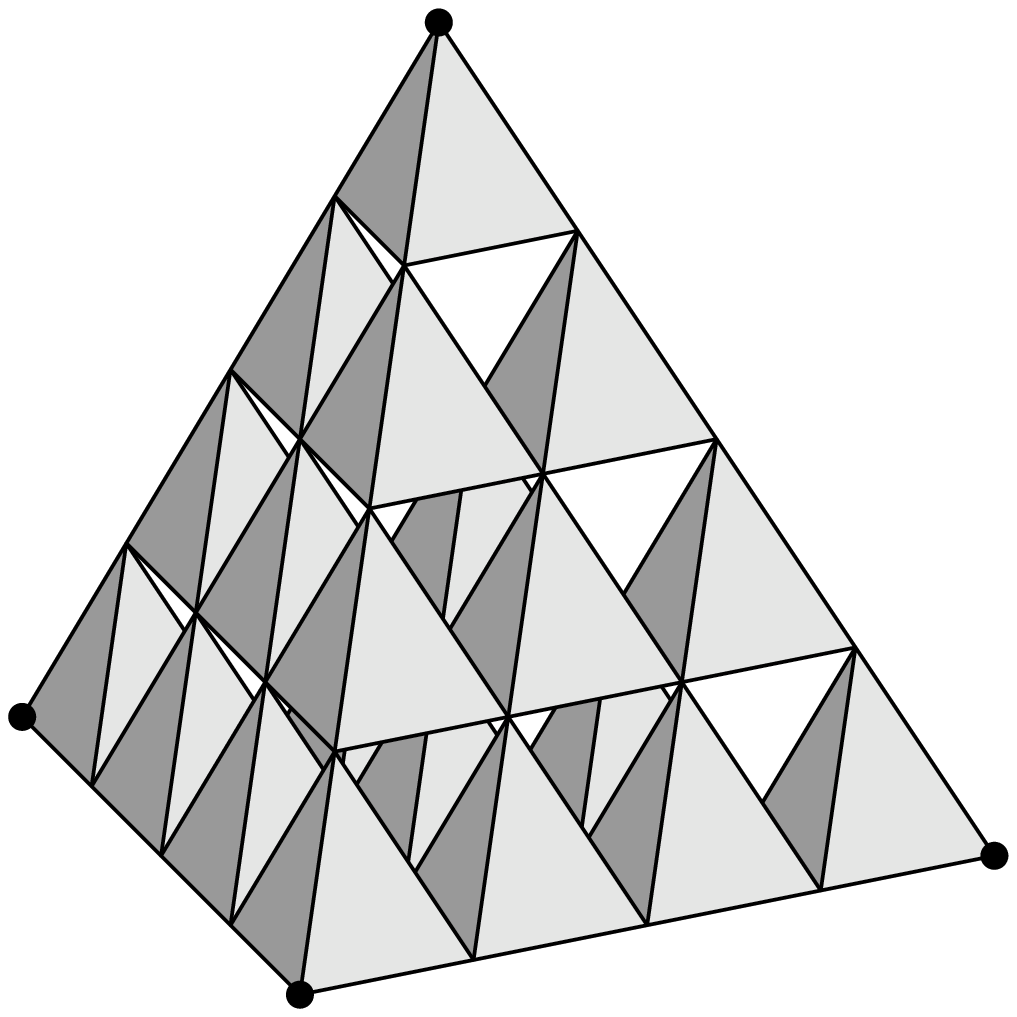}
\hspace{1cm}
\includegraphics*[width=5cm]{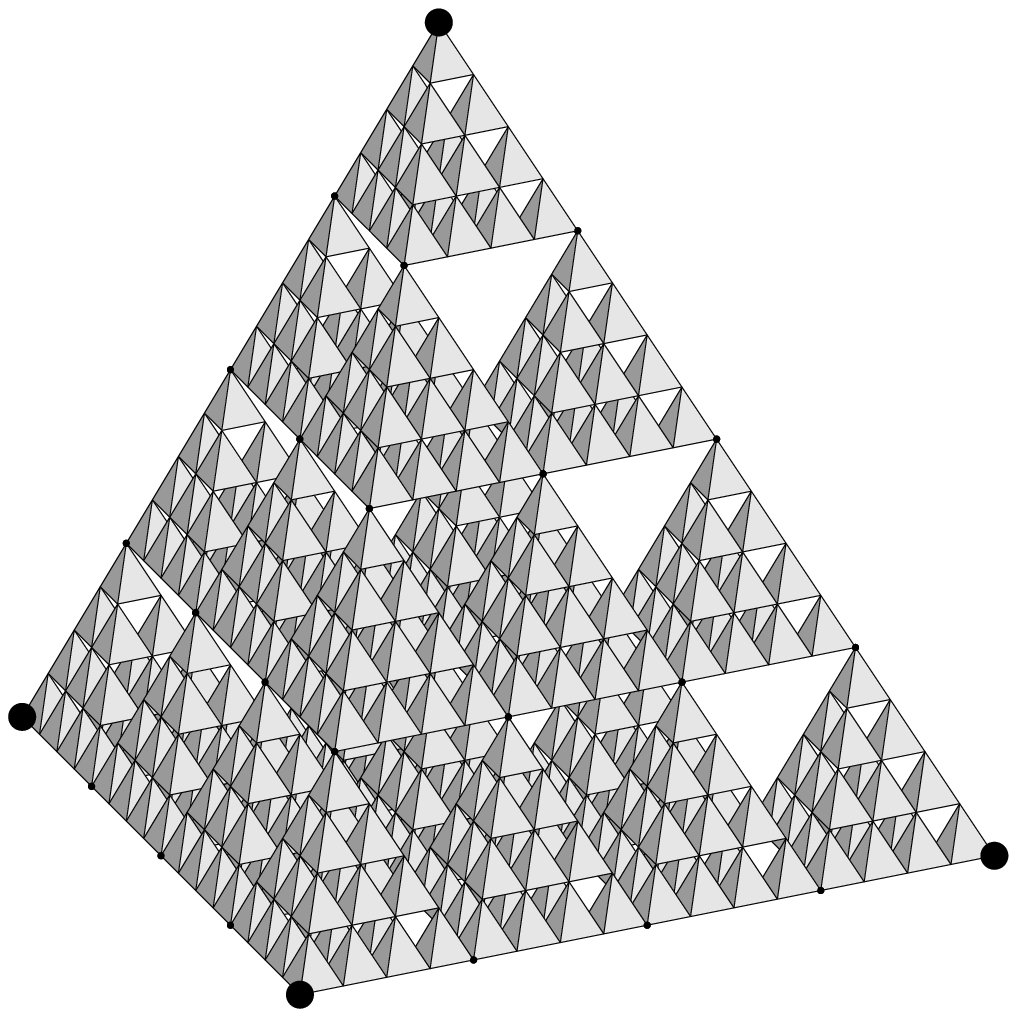}
\end{center}
\emph{Figure}~\ref{fig:34sierp}
\end{figure}
Here we have
\[ d(z)=\frac{z^4 (-486 + 9 z + 23 z^2 + 2 z^3)}
	{-104976 + 227448 z - 156168 z^2 + 31212 z^3 + 2958 z^4 - 887 z^5 - 41 z^6 + 2 z^7}  \]
and $\mu=20$, $\theta=4$, $\beta=1$, $\tau=\frac{4415}{113}$ and $\rho=\frac{883}{452}$. As
\[ d\bigl(\tfrac{9}{4}\bigr) = 1 \qquad\text{and}\qquad d(2) = \tfrac{36}{53} \]
the Julia set $\cj$ is a Cantor set.
If the constructed self-similar graph has an origin vertex
then the $1$-periodic function $\sigma$ is not constant.
\end{exmp}

\section*{Acknowledgements}

The authors want to thank Peter~Grabner, Klaus~Schmidt and Wolfgang~Woess
for fruitful discussions and financial support.
Especially, the results in Section~\ref{sec:geometric} would not have been possible
without the contributions of Peter~Grabner and Wolfgang~Woess.

\bibliographystyle{plain}

\begin{thebibliography}{10}

\bibitem{alexander82density}
S.~Alexander and R.~Orbach.
\newblock Density of states on fractals: fractons.
\newblock {\em J. Physique Lettres}, 43:L625--L631, 1982.

\bibitem{barlow98diffusions}
M.~T. Barlow.
\newblock Diffusions on fractals.
\newblock In {\em Lectures on probability theory and statistics (Saint-Flour,
  1995)}, pages 1--121. Springer, Berlin, 1998.

\bibitem{barlow88brownian}
M.~T. Barlow and E.~A. Perkins.
\newblock Brownian motion on the {S}ierpi\'nski gasket.
\newblock {\em Probab. Theory Related Fields}, 79(4):543--623, 1988.

\bibitem{beardon91iteration}
A.~F. Beardon.
\newblock {\em Iteration of rational functions}.
\newblock Springer-Verlag, New York, 1991.

\bibitem{arous00large}
G.~Ben~Arous and T.~Kumagai.
\newblock Large deviations of {B}rownian motion on the {S}ierpinski gasket.
\newblock {\em Stochastic Process. Appl.}, 85(2):225--235, 2000.

\bibitem{bremaud99markov}
P.~Br{\'e}maud.
\newblock {\em Markov chains}.
\newblock Springer-Verlag, New York, 1999.

\bibitem{bruijn79asymptotic}
N.~G. de~Bruijn.
\newblock An asymptotic problem on iterated functions.
\newblock {\em Nederl. Akad. Wetensch. Indag. Math.}, 41(2):105--110, 1979.

\bibitem{flajolet90singularity}
P.~Flajolet and A.~Odlyzko.
\newblock Singularity analysis of generating functions.
\newblock {\em SIAM J. Discrete Math.}, 3(2):216--240, 1990.

\bibitem{fukushima92spectral}
M.~Fukushima and T.~Shima.
\newblock On a spectral analysis for the {S}ierpi\'nski gasket.
\newblock {\em Potential Anal.}, 1(1):1--35, 1992.

\bibitem{goulden83combinatorial}
I.~P. Goulden and D.~M. Jackson.
\newblock {\em Combinatorial enumeration}.
\newblock John Wiley \&\ Sons, New York, 1983.

\bibitem{grabner97functional2}
P.~J. Grabner.
\newblock Functional iterations and stopping times for {B}rownian motion on the
  {S}ierpi\'nski gasket.
\newblock {\em Mathematika}, 44(2):374--400, 1997.

\bibitem{grabner97functional1}
P.~J. Grabner and W.~Woess.
\newblock Functional iterations and periodic oscillations for simple random
  walk on the {S}ierpi\'nski graph.
\newblock {\em Stochastic Process. Appl.}, 69(1):127--138, 1997.

\bibitem{inninger01rational}
C.~Inninger.
\newblock {\em Rational iteration}.
\newblock Universit\"atsverlag Rudolf Trauner, Linz, 2001.
\newblock Dissertation, University of Linz, 2001.

\bibitem{jones96transition}
O.~D. Jones.
\newblock Transition probabilities for the simple random walk on the
  {S}ierpi\'nski graph.
\newblock {\em Stochastic Process. Appl.}, 61(1):45--69, 1996.

\bibitem{kigami93harmonic}
J.~Kigami.
\newblock Harmonic calculus on p.c.f.\ self-similar sets.
\newblock {\em Trans. Amer. Math. Soc.}, 335(2):721--755, 1993.

\bibitem{kigami01analysis}
J.~Kigami.
\newblock {\em Analysis on fractals}.
\newblock Cambridge University Press, Cambridge, 2001.

\bibitem{kigami93weyl}
J.~Kigami and M.~L. Lapidus.
\newblock Weyl's problem for the spectral distribution of {L}aplacians on
  p.c.f.\ self-similar fractals.
\newblock {\em Comm. Math. Phys.}, 158(1):93--125, 1993.

\bibitem{kroen02green}
B.~Kr{\"o}n.
\newblock Green functions on self-similar graphs and bounds for the spectrum of
  the {L}aplacian.
\newblock to appear in Ann. Inst. Fourier \textbf{52} (2002), no. 6.

\bibitem{kroen02growth}
B.~Kr{\"o}n.
\newblock Growth of self-similar graphs.
\newblock preprint, 2002.

\bibitem{lindstroem90brownian}
T.~Lindstr{\o}m.
\newblock Brownian motion on nested fractals.
\newblock {\em Mem. Amer. Math. Soc.}, 83(420):iv+128, 1990.

\bibitem{malozemov95pure}
L.~Malozemov and A.~Teplyaev.
\newblock Pure point spectrum of the {L}aplacians on fractal graphs.
\newblock {\em J. Funct. Anal.}, 129(2):390--405, 1995.

\bibitem{malozemov01self}
L.~Malozemov and A.~Teplyaev.
\newblock Self-similarity, operators and dynamics.
\newblock preprint, 2001.

\bibitem{odlyzko82periodic}
A.~M. Odlyzko.
\newblock Periodic oscillations of coefficients of power series that satisfy
  functional equations.
\newblock {\em Adv. in Math.}, 44(2):180--205, 1982.

\bibitem{rammal84random}
R.~Rammal.
\newblock Random walk statistics on fractal structures.
\newblock {\em J. Statist. Phys.}, 36(5-6):547--560, 1984.

\bibitem{rammal83random}
R.~Rammal and Toulouse.
\newblock Random walks on fractal structures and percolation clusters.
\newblock {\em J. Physique Lettres}, 44:L13--L22, 1983.

\bibitem{sabot01spectral}
C.~Sabot.
\newblock Spectral properties of hierachical lattices and iteration of rational
  maps.
\newblock preprint, 2001.

\bibitem{teufl01average}
E.~Teufl.
\newblock The average displacement of the simple random walk on the
  {S}ierpi\'nski graph.
\newblock to appear in Combin. Probab. Comput.

\bibitem{woess00random}
W.~Woess.
\newblock {\em Random Walks on Infinite Graphs and Groups}.
\newblock Cambridge University Press, Cambridge, 2000.

\end{thebibliography}

\end{document}